\documentstyle{amsppt}
\nologo
%
%
\font\b=cmr10 scaled \magstep4
\def\bigzerou{\smash{\kern-25pt\lower1.7ex\hbox{\b 0}}}
\hsize=360pt
\vsize=500pt
\hbadness=5000
\tolerance=1000
\NoRunningHeads

\def\z{\zeta}
\def\Q{\Bbb Q}

\def\ce{\hat{\Cal E}_W}
\def\cg{\widehat{\Bbb G_m}_W}

\def\al{\alpha}
\def\be{\beta}
\def\f{\frak}
\def\Fp{\Bbb F_p}
\def\F5{\Bbb F_5}

\def\Z{\Bbb Z}

\def\ep{\epsilon}

\def\t{\theta}
\def\no{\noindent}
\def\la{\lambda}
\def\p{\Phi}
\def\vp{\varphi}
\def\i{\infty}
\def\ti{\times}
\def\si{\sigma}

\def\Ga{\Gamma}
\def\ga{\gamma}

\def\op{\oplus}
\def\ot{\otimes}
\def\part{\partial}
\def\va{\varprojlim}
\magnification=\magstep1
\topmatter
\title Elliptic units in $K_2$  \endtitle
\author
Kenichiro Kimura
\endauthor
\footnote{The author was partially supported by JSPS Postdoctoral Fellowships
for Research Abroad from 2000-2001.}

\affil
Institute of Mathematics, University of Tsukuba,
Tsukuba, Ibaraki, 305 Japan
email: kimurak\@math.tsukuba.ac.jp
\endaffil

\endtopmatter
\vskip -3ex
\noindent\hskip 2em
\vskip 3ex
\CenteredTagsOnSplits
%
\document

\baselineskip=13pt

abstract: We present certain norm-compatible systems in $K_2$ of 
function fields of some CM elliptic curves. We demonstrate that 
these systems have some properties similar to elliptic units.


\subhead{\S1 Introduction}\endsubhead
Elliptic units are systems of units
in abelian extensions of imaginary quadratic fields which are
defined by evaluating some functions on CM elliptic curves
at torsion points. These units play an essential role in Iwasawa
theory of imaginary quadratic fields (\cite{dS},\cite{Ru} etc).
In this paper we present certain systems of $K_2$ of function fields
of CM elliptic curves. These systems are considered by Scholl
(\cite{Sch}) in a more general setting. We show that these
systems have two similar properties to elliptic units.

\no First,  in $\S2$ we demonstrate that we can apply the theory 
of Coleman power series for $K_2$ of 2--dimensional local fields 
constructed by Fukaya (\cite{Fu}) to our system.
We show that the ``Coleman power series'' of our systems, which are
elements in $K_2$ of some power series rings, are explicitly
described.   In $\S3$ we modify our systems so that they are in 
the kernel of tame symbol maps. Using the result of $\S2$ we will 
give a sufficient condition for the images of those systems to be
non-zero under Chern class maps to Galois cohomology groups.

\no Second, in the appendix we show that our systems form
Euler systems similar to those of elliptic units.

By \cite{Sch} one can see that systems like ours also exist
in higher $K$-groups of function fields of products of elliptic curves.
This fact seems to suggest (to this author) that there is a further 
generalization of the theory of elliptic units and Coleman power series
to higher $K$-groups.

\no The following question is asked by the referee:

\no {\it Is there interpretation of the function which appears in corollary
3.4 in terms of L-function?}

\subsubhead{Acknowledgment}\endsubsubhead  The author would like
to express his gratitude to Prof. K. Kato for valuable suggestions
and to Dr.T. Fukaya for showing her paper at an early stage. He is also 
grateful to 
the referee for valuable comments on the first version of this paper.
 This work
was carried out while the author stayed at Department of Mathematics,
the University of Chicago. He is thankful to the people at Chicago 
for their hospitality.

\subhead{Notations}\endsubhead

 $K$: an imaginary quadratic field, $\quad$ 4

 $E$: a CM elliptic curve, $\quad$ 4

 $F$: an abelian extension of $K$, $\quad$ 4

 $\psi$:  the Hecke character of $F$ associated to $E$, $\quad$ 5

 $\phi$: a Hecke character of $K$ such that $\psi=\phi\circ N_{F/K}$, $\quad$ 5

 $\frak f$: conductor of $E$, $\quad$ 5

 $p$: a rational prime, $\quad$ 5

 $\frak p$: a prime factor of $p$ in $\Cal O_K$, $\quad$ 5

 $\pi$: $\phi(\frak p)$, $\quad$ 5

 $a$: an integer relatively prime to $6p\frak f$, $\quad$ 5

 $g_a$: a rational function on $E$, $\quad$ 5

 $v_n$: a generator of $E[\pi^n]$, $\quad$ 5

 $T_y$: translation by a point $y$, $\quad$ 5

 $F_n$: $F(E[\pi^n])$, $\quad$ 5

 $\al_n$: an element of $K_2(F_n(E))$, $\quad$ 5

 $W$: the Witt ring of $\overline{\Bbb F_p}$, $\quad$ 6

 $\frak P$: a prime of $F$ over $p$, $\quad$ 6

 $\hat{\Cal E}_W$: the formal completion of  a projective smooth model of $E$, $\quad$ 6

 $\theta$: an isomorphism $\cg \to \ce$, $\quad$ 6

 $\Phi$: Frobenius of $W$, $\quad$ 6

 $t$: a parameter on $\ce$, $\quad$ 6

 $s$: the standard parameter on $\cg$, $\quad$ 6

 $(\z_n)_{n\leq 1}$: a basis of $\Z_p(1)$, $\quad$ 6

 $H$: $p$-adic completion of $W((s))$, $\quad$ 6

 $H_n$: $H(s_n,\,\z_n)$, $\quad$ 6

 $S$: $\Cal O_H[[\ep -1]]$, $\quad$ 6

 $S'$: $S[(\ep-1)^{-1}]$, $\quad$ 6

 $\nu,\,\,\sigma,\,\,\varphi$: endomorphisms of $S$, $\quad$ 6

 $h$: an isomorphism ${K_2(S' )^{\wedge}}^{Nr_\vp=1}
\simeq \varprojlim_n K_2(H_n)^{\wedge}$, $\quad$ 7

 $G_k$: $Gal(k^{sep}/k)$ for a field $k$, $\quad$ 7

 $H^i(k,\,M)$: $H^i(G_k,\,M)$ for a $G_k$ module $M$, $\quad$ 7

 $\lambda, \part$: homomorphisms $\varprojlim_nK_2(H_n)^{\wedge} \to \Omega^1_{H_1}$, $\quad$ 7

 $q$ : $s+1$, $\quad$ 8

 $\be_n$: ${\t^{\Phi^{-n}}}^* \al_n$, $\quad$ 8

 $D$: $q\frac{d}{ds}$, $\quad$ 9

 $F_n'$: $F(E[a\pi^n])$, $\quad$ 9

 $\al_n^\gamma$: $T_{\gamma}^*\al_n$, $\quad$ 9

 $\al_n'$: an element of $K_2(F_n(E))$, $\quad$ 9

 $\beta_n'$: ${\t^{\p^{-n+1}}}^*\al'_n$, $\quad$ 10

 $g_a^{\ga}$: $T_{\ga}^*g_a$, $\quad$ 11

\subhead{\S 2 norm-compatible systems in $K_2$}\endsubhead

A general reference for this section is Chapter II of \cite{dS}.
We fix an imaginary quadratic field $K$.  Let $O_K$ denote the
ring of integers of $K$. Let F be an abelian extension of $K$ and 
let $E$ be an elliptic curve over $F$ which satisfies
the following conditions:

(i)  $E$ has complex multiplication by $O_K$, and

(ii) $F(E_{tors})$ is an abelian extension of $K$.

$\quad$

\no Here $F(E_{tors})$ denotes the field generated by coordinates of 
torsion points of a model of $E$ over $F$.
If $\psi$ is the Hecke character of $F$ associated to $E$,
then the condition (ii) is equivalent to the existence of a Hecke
character $\phi$ of $K$ such that $\psi=\phi\circ N_{F/K}.$

\no Let $\f f$ be the conductor of $E$ and let $p>3$ be a rational prime
which is relatively prime to $\f f$ and which splits completely
in $F$. Let $p=\f p\bar{\f p}$ be the decomposition
of $p$ in $O_K$. We write $\pi=\phi(\f p)$. Note that $\f p=(\pi).$
We assume that $Gal(F(E[p^n]))/F\simeq (O_K/(p^n))^*$ for all $n\geq 1$.

\no Let $a$ be an integer such that $(a,\,\,6p\f f)=1$.
Then there is a function $g_a\in F(E)$, which is called the Kato-Siegel
function in \cite{Sch}. We now state the properties of $g_a$ which
will be used in the subsequent discussions.

(a)  div$g_a=a^2(0)-E[a]$.

(b) For any isogeny $j:\,\,E\to E$ defined over $\bar{F}$ such that
$(\text{deg}\, j,\,a)=1$, $j_*g_a=g_a$.

(c) $[a]_*g_a=1$.

(d) $[-1]^*g_a=g_a$. (This is a special case of (b)).
$\quad$

\no Here for an element $u$ of $O_K$, $[u]$ stands for an element
of End($E$), which is justified by the condition (i). 
In our situation, the property (b) is equivalent to the following
distribution relation
$$\prod_{j(y)=x}g_a(y)=g_a(x).$$
We fix a system $(v_n)_{n\ge 1}$ of torsion points of $E(\bar{F})$ such that
for each $n\geq 2$, $v_n$ is a generator of $E[\pi^n]$ and that
$[\pi]v_n=v_{n-1}$.

\no For a point $y$ on $E$ we denote by $T_y:\,\,E\to E$
the translation by $y$. For $n\geq 1$ consider the element

$$\al_n:=\{g_a,\,T_{v_n}^*g_a\}\in K_2(F_n(E)).$$
Here $F_n:=F(E[\pi^n])$ and $F_n(E)$ is the function field of
$E_{F_n}=E\ti_{F}F_n$. This element is considered by Scholl(\cite{Sch},
Appendix) in a more general setting.
The following lemma claims that the system $(\al_n)$ forms
an norm-compatible system.

\proclaim{Lemma 2.1} For $n\geq 2$, the equality
$$[\pi]_*N_{F_n(E)/F_{n-1}(E)}\al_n=\al_{n-1}$$
holds. \endproclaim

\demo{Proof} This is another way of describing (EN1) in loc.cit.
We will give a proof by calculation of symbols.
$$N_{F_n(E)/F_{n-1}(E)}\al_n=\{g_a,\,\prod_{v\in E[\pi]}
T_{v_n+v}^*g_a\}=\{g_a,\,[\pi]^*(T_{v_{n-1}}^*g_a)\}$$ where
the second equality holds by the distribution relation mentioned above, 
and
$$[\pi]_*\{g_a,\,[\pi]^*(T_{v_{n-1}}^*g_a)\}=
\{[\pi]_*g_a,\,T_{v_{n-1}}^*g_a\}=\{g_a,\,T_{v_{n-1}}^*g_a\}.$$
\enddemo

We will apply the theory of Coleman power series for $K_2$
constructed by Fukaya \cite{Fu} to this system.  As at present her
theory treats only the multiplicative group, we need to construct 
an isomorphism from $E$ to the multiplicative group.  For this
we will take the formal completion of $E$ along the $0$-section.

\no We will follow the argument of \S 4.3,Chapter II of \cite{dS} rather
closely. Let $F'=F(E[\bar{\pi}^\i])$ and let $R$ be the completion of
$O_{F'}\ot_{O_K}{O_K}_{\f p}$. Since the prime $\f p$
finitely decomposes in $F'$, $R$ is a direct sum of local
rings. Fix a direct summand of $R$ and denote it by $W$. 
Let $\f P$ be the prime of $F$ over $\f p$ which corresponds to $W$
and let $\Cal E$ be a smooth projective model
of $E$ over ${O_F}_{\f P}=\Z_p$. Let $\ce$ be the formal
completion of $\Cal E\ti_{\Z_p} W$ along the $0$-section. Then 
by Proposition 4.3 in Chap.2 of \cite{dS}, there is an isomorphism

$$\t:\,\,\cg \to \ce$$
which satisfies the equality $[\pi]\circ \t=\t^\p\circ [p]$.
Here $\p$ is the Frobenius automorphism of $W$. We fix a parameter
$t$ of $\ce$ and let $s$ be the standard parameter of  $\cg$.

\no For $n\ge 1$, define $\z_n\in \overline{\text{Frac}(W)}$
by the equality $t(v_n)=\t^{\p^{-n}}(\z_n-1)$. Then we see that 
the system $(\z_n)_{n\geq 1}$ is a basis of $\Z_p(1).$

Now we will review the theory of Fukaya in the form we will use in
the subsequent discussions.

\no Let $H$ be the $p$-adic completion of $W((s))$, the function
field of $\cg$.
Consider a tower of field extensions $(H_n)_{n\geq 1}$
where $H_n:=H(s_n,\,\z_n)$ and $(s_{n-1}+1)=(s_n+1)^p$ for $n\geq 1$. 
$s_0=s$ by convention.

\no We set $S:=O_H[[\epsilon-1]]$ to be the ring of formal power
series and write $S'=S[(\ep-1)^{-1}]$. Consider the ring homomorphisms
$\nu,\,\,\si$ and $\vp:\,\,S'\to S'$ defined as follows:

\no $\nu$ is defined by $\nu(\ep)=\ep^p$ and $\nu|_{O_H}=$
the identity. The homomorphism $\si$ is defined by
$\si|_{O_H}=[p]^*\ot \p$ and $\si(\ep)=\ep.$
The homomorphism $\vp$ is defined by the composition: $\vp=\nu\circ \si$.
The following theorem is a part of Theorem 1.5 of \cite{Fu}.
It asserts the existence of Coleman power
series:

\proclaim{Theorem 2.2} There is a canonical isomorphism

$$h:\,\,{K_2(S' )^{\wedge}}^{Nr_\vp=1}
\simeq \varprojlim_n K_2(H_n)^{\wedge}$$
induced by

$$h_n:\,\,S'\to H_n,\,f(\ep)\mapsto (\si^{-n}f)(\z_n).$$
Here $Nr_\vp$ is the norm map on $K_2(S')$ induced by the finite
morphism $\vp$ and $*^{Nr_\vp=1}$ means the fixed part of $*$ under 
$Nr_\vp$. The map $\si^{-n}:\,\,H\hookrightarrow H_n$ is the composition
of $id\ot\p^{-n}:\,\,H\to H$ and the inclusion given by 
$s\mapsto s_n$.
 The symbol $\wedge$ means certain completion.
\endproclaim

\no Also there is a fact which corresponds to the explicit reciprocity
law. It asserts the coincidence with $-$ sign of the two homomorphisms 
$\la$ and $\part$ defined as follows. 

\no For a field $k$,  we write $G_k$ for $Gal(k^{sep}/k)$. For a $G_k$
module $M$, we write $H^i(k,\, M)$ for $H^i(G_k,\,M).$


$$\align\la:\,\,\varprojlim_nK_2(H_n)^{\wedge}&\overset{cl}\to\to
\va_nH^2(H_n,\,(\Z/p^n\Z)(2))\\
&\overset{\cup(\z_n)^{-1}}\to\longrightarrow
\va_nH^2(H_n,\,(\Z/p^n\Z)(1))\\
&\overset{Cor_{H_n/H_1}}\to\longrightarrow\va_nH^2(H_1,\,(\Z/p^n\Z)(1))\\
&=H^2(H_1,\,\Z_p(1))\\
&\to H^2(H_1,\,\Bbb C_p(1))\simeq \Omega^1_{H_1}.\endalign$$
The last isomorphism is due to Hyodo \cite{Hy}.
Here $ \Omega^1_{H_1}$ is defined to be
$\Q\ot \Omega_{O_{H_1}}$ where
$\Omega_{O_L}^1:=\underset{n}\to\varprojlim \Omega^1_{O_{H_1}/\Z}
/p^n\Omega^1_{O_{H_1}/\Z}$.

The homomorphism $\part$ is defined using the map $h$ in Theorem 2.2:

$$\align
\part:\,\,\va_n K_2(H_n)^{\wedge}
&\overset{h^{-1}}\to{\longrightarrow} K_2(S')^{Nr_\vp=1}\\
&\overset{d\log}\to{\to}\Omega_S^2(\log)
=\frac{S}{\ep-1}d\log(q)\wedge d\log(\ep)\\
&\overset{(\wedge d\log(\ep))^{-1}}\to{\longrightarrow}
\frac{S}{\ep-1}d\log(q)\\
&\to \Omega^1_{H_1}\endalign$$
where $q=s+1$ and the last map is defined by
$$f(\ep)d\log(q)\mapsto \frac 1p(\si^{-1}f)(\z_1)d\log(q).$$
The modules $\Omega^1_S(\log)$ and $\Omega^2_S(\log)$ are
defined as follows:

\no Let
$$\Omega^1_{S/\Z}(\log):=(\Omega^1_{S/\Z}\op S\ot_\Z {S'}^*)/N$$
where $N$ is the $S$ submodule generated by elements
$(-da,\, a\ot a)$ for $a\in S\cap {S'}^*$. The class $(0,\,1\ot a)$
is denoted by $d\log (a).$ The module $\Omega^2_S(\log)$ is defined
by
$$\Omega^2_S(\log)=\va_n\wedge^2
\Omega^1_{S/\Z}(\log)/p^n\wedge^2\Omega^1_{S/\Z}(\log).$$ 
By proposition (19.3.5) of \cite{Gr}, $S$ is formally smooth over
$\Bbb Z$, so by corollary (20.4.11) of loc.cit. 
$\Omega^1_{S/\Z}(\log)/p^n\Omega^1_{S/\Z}(\log)$ is a free $S/p^nS$ module
on $d\log (\ep-1)$ and $d\log(q)$. Hence $\Omega^2_S(\log)$ is
a free $S$ module on $d\log(q)\wedge d\log(\ep-1).$

\no  The map $d\log:\,\,K_2(S')\to \Omega_S^2(\log)$ is described as 
$\{f,g\}\mapsto d\log(f)\wedge d\log(g)$.

\proclaim{Theorem 2.3 (Theorem 2.4 in \cite{Fu})}
$$\la=-\part.$$
 \endproclaim

 For $n\ge 1$, set
$\be_n:={\t^{\p^{-n}}}^*\al_n
\in K_2(H_n).$ Since $[p]_*\t^*={\t^\Phi}^*[\pi]_*$ 
it follows that $(\beta_n)_n$ is an element
of $ \varprojlim_n K_2(H_n)$ where projective limit is taken
with respect to norm maps.

\no The following is our main result:

\proclaim{Theorem 2.4} 
We have the equality
$$h^{-1}((\beta_n)_n)=B:=\{\t^*g_a(s),\,\t^*g_a(s[+]_{\cg}\ep-1)\}.$$
\endproclaim

\demo{Proof} First of all, note that since $g_a(t)\in t^{1-a^2}W[[t]]^*$,
$\t^*g_a(s)\in s^{1-a^2}W[[s]]^*\subset O_H^*$. 
Since $\t^*g_a(s[+]_{\cg}\ep-1)$ is of the form
$\t^*g_a(s)+(\ep-1)f(s,\ep-1)$ with $f(s,\ep-1)\in S$, it belongs to 
$S^*$. It can also be
checked directly that $h_n(B)=\beta_n$ for $n\geq 1.$ What remains is
to show that $Nr_{\vp}B=B.$

$$\align Nr_{\vp}B=&Nr_{\si}\circ Nr_{\nu}B\\
    =&Nr_{\si}\circ \{\t^*g_a,\,\prod_{\z^p=1}\t^*g_a(s[+]_{\cg}\z\ep-1)\}\\
    =&Nr_{\si}\circ \{\t^*g_a,\,\prod_{\z^p=1}\t^*g_a(s[+]_{\cg}\ep-1
   [+]_{\cg}\z-1) \}\\
   =&Nr_{\si}\circ \{\t^*g_a,\,\prod_{\z^p=1}g_a(\t(s[+]_{\cg}\ep-1)
   [+]_{\ce}\t(\z-1)) \}\\
   =&Nr_{\si}\circ \{\t^*g_a,\,\prod_{v\in E[\pi]}g_a(\t(s[+]_{\cg}\ep-1)
   [+]_{\ce}v) \}\\
   =&Nr_{\si}\circ \{\t^*g_a,\,g_a([\pi](\t(s[+]_{\cg}\ep-1))) \}\\
   =&Nr_{\si}\circ \{\t^*g_a,\,g_a(\t^\p([p](s[+]_{\cg}\ep-1))) \}\\
=&Nr_{\si}\circ \{\t^*g_a,\,[p]^*({\t^\p}^*g_a)(s[+]_{\cg}\ep-1)
\}.\endalign$$
Since $\si=[p]^*\ot \p$, $Nr_{\si}=\p^{-1}\circ [p]_*.$
From this and the equality
$[p]_*(\t^*g_a)={\t^\p}^*[\pi]_*g_a={\t^\p}^*g_a$ it follows that
$Nr_\vp B=B.$ \qed \enddemo

\no Let $D$ be a differential operator
 on $H$ defined by $D=q\frac{d{}}{ds}$. Then by straightforward calculation
we obtain the following corollary.
\proclaim{Corollary 2.5}
$$\part((\beta_n)_n)=
\frac 1p\p^{-1}(D\log(\t^* g_a)(s)\ti
D\log(\t^* g_a)(s[+]_{{\cg}} \z_1-1))d\log q.$$
\endproclaim

 \subhead \S3. Elements in the kernel of tame symbol map\endsubhead

For $n\ge 1,$ let $F'_n:=F(E[a\pi^n])$. For a point $\ga\in E[a]$,
let $\al_n^\ga:=T^*_{\ga}\al_n\in K_2(F'_n(E))$.
We set

$$\al_n':= N_{F'_n/F_n}[\pi]_*\underset{\ga\in E[a]}\to\sum \al_n^\ga.$$
We can check that for $n\geq 2$,
$$ [\pi]_*N_{F_n(E)/F_{n-1}(E)}\al_n'=\al_{n-1}'$$ 
in a similar way as in the case of $\al_n$.
We have the following exact sequence (\cite{Qu}): 
 $$\underset{x\in E(\bar{F})}\to\op K_2(\kappa_x) \to
K_2(E_{F_n})\overset{r}\to \to K_2(F_n(E))\overset
{\underset{x\in E(\bar{F})}\to\op t_x}\to{\longrightarrow}
\underset{x\in E(\bar{F})}\to\op \kappa_x^*
\overset{N}\to\to F_n^*$$
where $r$ is the restriction, $N$ is the norm and for $x\in E(\bar{F})$,
$t_x$ is defined by $$\{f,\,g\}\mapsto
(-1)^{ord_xford_xg}\frac{f^{ord_xg}}{g^{ord_xf}}(x).$$

\no By this exact sequence we see that
 $[\pi^n]^*[\pi^n]_*\underset{\ga\in E[a]}\to\sum \al_n^\ga
=\underset{u\in E[a\pi^n]}\to\sum T^*_u \al_n$ is
in the kernel of$\underset{x\in E(\bar{F})}\to\op t_x$, and since
  the map $[\pi^n]:\,\,E\to E$ is etale, it follows that
  $\underset{x\in E(\bar{F})}\to\op t_x([\pi^{n-1}]_*\al_n')=0.$

\no Hence there is an element $\Cal A\in K_2(E_{F_n})$ such that $r(\Cal A)
=[\pi^{n-1}]_*\al_n'.$ Since $K_2$ of number fields is torsion (\cite{Ga}), 
$\Cal A$ is well defined up torsion.

 \proclaim{Lemma 3.1}Let $V_p(E)(1):=T_p(E)(1)\ot \Q$ and set
 $F_\i:=F(E[\pi^\i]).$ Then
 $$(V_p(E)(1)/T_p(E)(1))^{G_{F_\i}}=0.$$
 \endproclaim
 \demo{Proof} There is a decomposition $T_p(E)(1)=
 T_{\f p}(E)(1)\oplus T_{\bar{\f p}}(E)(1)$, so we will treat
 each factor separately. The factor $ T_{\f p}(E)(1)$ is isomorphic
 to $\Z_p(1)$ as a  $G_{F_\i}$ module. By corollary 1.7
 in Chap.II of \cite{dS}, the fields $F_\i$ and $F(E[\bar{\pi}^\i])$
 are linearly
 disjoint over $F$. Since $\Z_p(1)\simeq
 T_{\f p}(E)\ot T_{\bar{\f p}}(E)$, the fields
 $F_\i$ and the field $F(\z_{p^\i})$, generated by $p$-power
roots of unity, are also linearly disjoint over $F$.
 Hence the assertion holds for this factor.

 \no As for $T_{\bar{\f p}}(E)(1)$, it is isomorphic to
 $\Z_p(2)$ as a $G_{F_\i}$ module. Since $p>3$,
 the assertion holds for this factor as well. \qed\enddemo

 \no By the lemma, $H^1(F_n,\,T_p(E)(1))$ has no torsion
 and we have a well defined class
 $\tilde{cl}([\pi^{n-1}]_*\al'_n):=cl(\Cal A)\in H^1(F_n,\,T_p(E)(1)).$

 We set $\beta_n':={\t^{\p^{-n+1}}}^*\al'_n \in K_2(H_n).$
 Since all the points in $E[a]$ are rational over Frac$W$ 
and ${\t^\Phi}^*[\pi]_*=[p]_*\t^*$, 
 we have the equality $$\beta_n'=[F_0':F][p]_*
 \underset{\ga\in E[a]}\to\sum {\t^{\p^{-n}}}^*\al_n^\ga.$$
 One can check that $(\beta'_n)_n$ form an element of
 $\varprojlim_n K_2(H_n)$ in a similar way to the case of $\be_n$.

 \no It can also be shown that $\t^*([\pi^{n-1}]_*\al_n')=
 [p^{n-1}]_*\be_n'$. By  (5.10) in \cite{Fu}, we have a commutative
 diagram
 $$\CD
  \va_nK_2(H_n)^{\wedge} @>cl>> \va_nH^2(H_n,\,\Z/p^n\Z(2))
  @> C_1>>H^2(H_1,\,\Z_p(1))\\
  @V[p^{n-1}]_*VV @V\text{Cor}_{H_n/H_1(\z_n)}VV @V\text{id}VV\\
  \va_nK_2(H_1(\z_n))^{\wedge} @>cl>> \va_nH^2(H_1(\z_n),\,\Z/p^n\Z(2))
  @>C_2 >>H^2(H_1,\,\Z_p(1))
  \endCD $$
Here $C_1=\text{Cor}_{H_n/H_1}\circ\cup(\z_n)^{-1}$ and 
$C_2=\text{Cor}_{H_1(\z_n)/H_1}\circ \cup(\z_n)^{-1}$.
By this diagram it follows that if $\part((\be_n')_n)\neq 0$, then
$cl([p^{n-1}]_*\be_n')\neq 0$ for sufficiently large $n$. 

\no It implies that $0\neq cl([\pi^{n-1}]_*\al'_n)\in H^2({F_n}_\f P(E),\,\Z/p^n\Z(2))$
where $ {F_n}_\f P$ is completion of $ {F_n}$ by the prime over $\f P.$

Consider the restriction of $\tilde{cl}([\pi^{n-1}]_*\al_n')$ to
$H^1({F_n}_{\f P},\,T_p(E)(1))$.
We need a splitting
$H^1({F_n}_{\f P},\,T_p(E)(1))\to  H^2(E_{{F_n}_{\f P}},\,\Z/p^n\Z(2))$
and the following lemma enables us to define such a splitting.

\proclaim{Lemma 3.2} $\sharp H^2({F_n}_{\f P},\,\Z/p^n\Z(2))
=\sharp H^0({F_n}_{\f P},\,\Z/p^n\Z(-1))$ is bounded
with respect to $n$.\endproclaim

\demo{Proof} Considering the Weil pairing $E[\pi^n]\ti
E[\bar{\pi}^n]\to \mu_{p^n}$, we see that if
${F_{n}}_{\f P}\supset F_{\f P}(\mu_{p^m})$ for a number $m\leq n$,
 then
$ {F_{n}}_{\f P}\supset
F_{\f P}(E[\bar{\pi}^m]).$
Since $\f P$ is finitely decomposed
in $F(E[\bar{\pi}^\i])$, the residue field of
$F_{\f P}(E[\bar{\pi}^\i])$is $\overline{\Fp}$. As ${F_n}_{\f P}$
is totally ramified over $F_{\f P}$, it follows that $m$ is bounded.
\qed\enddemo

\no By this lemma we have a number $M$ which is a power of $p$ and
a map
$$u_n:\,\,
H^1({F_n}_{\f P},\,T_p(E)(1))\to  H^2(E_{{F_n}_{\f P}},\,\Z/p^n\Z(2))$$
for each $n\ge 1$ such that 

$$ u_n(\tilde{cl}([\pi^{n-1}]_*\al_n') =Mcl([\pi^{n-1}]_*\al_n').$$

\no Noting that the action of the Frobenius $\p$ on $E[a]$ is by $\pi$,
 one can show the following in a similar way to the proof of Theorem 2.4.
 \proclaim{Theorem 3.3} For $\ga\in E[a]$, set $g_a^\ga:=T_\ga^*g_a.$
 Then we have the equality
 $$h^{-1}((\be_n'))=B':=[F_0':F]
 [p]_*\sum_{\ga\in E[a]}\{\t^*g_a^\ga,\,
 \t^*g_a^\ga(s[+]_{\cg} \ep-1)\}.$$\endproclaim

\no Since $d\log$ and the pull back $[p]^*$ commutes on $S'$, it also
 commutes on $K_2(S')$. As
 $$\align[p]^*B'=&[F_0':F]\sum_{\z^p=1}\sum_{\ga\in E[a]}\\
 &\{\t^*g_a^\ga(s[+]_{\cg}\z-1),
 \t^*g_a^\ga(s[+]_{\cg}\z-1[+]_{\cg} \ep-1)\}\endalign$$
 we see that
 $$\align &[p]^*\part((\be_n'))\\
 =&\frac{[F_0':F]}{p}\p^{-1}(\sum_{\z^p=1}\sum_{\ga\in E[a]}
 D\log \t^*g_a^\ga(s[+]_{\cg}\z-1)\\
 &\ti D\log \t^*g_a^\ga(s[+]_{\cg}\z-1[+]_{\cg} \z_1-1))d\log q.\endalign$$
 We write the last term as $\frac{[F_0':F]}{p}\p^{-1}G(s)d\log q.$ 
 
 Let $\la:\,\,\widehat{E}\to \widehat{\Bbb G}_a$ be a logarithm defined
 over $F$ and let
 $z$ be the parameter of $\widehat{\Bbb G}_a$ such that
 $z=\la(t)=t+\cdots $. We have the equality 

 $$\Omega_p\frac{d{}}{dz}\circ \la_*\t_* =\la_*\t_*\circ D$$
 where $\Omega_p\in W^{*}$ is such that $t=\t(s)=\Omega_ps+\cdots$.
 Let $\xi:\,E(\Bbb C)\to \Bbb C/L$ be the inverse of the analytic
 uniformization associated to $dz$. Here $L$ is the period lattice
 of $dz$. For a point $v\in E[\pi]$, $\frac{d{}}{dz}\la_*g_a^{\ga}
 (z[+]_{\widehat{\Bbb G}_a}\la(v))$ is the Taylor expansion
 of $\frac{d{}}{dz}\xi_*g_a^{\ga}
 (z+\xi(v))$ at $z=0.$ Hence by looking at 
 $G(s)$ we obtain the following corollary.
 
 \proclaim{Corollary 3.4}
 
 If the function
 $$\sum_{v\in E[\pi]}\sum_{\ga\in E[a]}
 \frac{d{}}{dz}\log (\xi_*g_a^\ga)(z+\xi(v))
\ti \frac{d{}}{dz}\log (\xi_*g_a^\ga)(z+\xi(v+v_1))$$
on $\Bbb C/L$ is nonzero, then $\tilde{cl}([\pi^{n-1}]_*\al'_n)$
is nonzero for sufficiently large $n$. 
 \endproclaim

\subhead{Appendix: Euler Systems}\endsubhead
Let $\Cal S$ be the set of ideals of $O_K$ which are divisible
only by primes which split completely in $F$ and which are
relatively prime to $\f f a$. For each $\f m\in \Cal S$, Let
$F_{\f m}:=F(E[\f m]).$ Fix a system of torsion points
$(v_{\f m})_{\f m\in \Cal S}$ such that

(1) $ v_{\f m}$ is a generator of $E[\f m]$ for each $\f m\in \Cal S.$

(2) If $\f m_2=\f m_1\f n$ with $\f m_2\in \Cal S$, then
$v_{\f m_1}=\phi(\f n)v_{\f m_2}$.

\no We set $\al_{\f m}:=\{g_a,\,T_{v_{\f m}}^*g_a\}\in K_2(F_{\f m}(E))$
and set $\al_{\f m}':=[\phi(\f m)]_*\underset{\ga\in E[a]}\to\sum
T_\ga^*\al_{\f m}.$
Then it is seen that $\al_{\f m}'\in \Ga(E_{F_{\f m}},\,\Cal K_2)$
 in the similar way to  the argument in the beginning
 of \S 3.  So that there is a class
$\tilde{cl}(\al_{\f m}')\in H^1(F_{\f m},\,T_p(E)(1)).$
In the similar way to the proof of (EN1) in (\cite{Sch}, Appendix)
one obtains the following result:

\proclaim{ Proposition A.1} Let $\f m\in \Cal S$ and assume that
$\f l|\f m$ for a prime $\f l$. Then the following equality holds:
$$N_{F_{\f m}/F_{\f m/\f l}}\al_{\f m}'
        =\cases &\al_{\f m/\f l}'\quad \f l^2|\f m\\
                 &(1-[\phi(\f l)]_*Fr_{\f l}^{-1})\al_{\f m/\f l}'
                 \quad \f l^2\nmid \f m.\endcases$$
\endproclaim

\no The classes $\tilde{cl}(\al_{\f m}')$ satisfy
the same relation as $\al_{\f m}'$. This implies the classes 
$\tilde{cl}(\al_{\f m}')$ form an Euler system.

\enddocument